\documentclass{article}
\usepackage[english]{babel}

\newcommand{\vertiii}[1]{{\left\vert\kern-0.25ex\left\vert\kern-0.25ex\left\vert #1 
		\right\vert\kern-0.25ex\right\vert\kern-0.25ex\right\vert}}

\newcommand{\vvertiii}[1]{{\vert\kern-0.25ex\vert\kern-0.25ex\vert #1 
		\vert\kern-0.25ex\vert\kern-0.25ex\vert}}

\newcommand{\dproof}{\noindent {Proof.} \quad}
\newcommand{\fproof}{\hfill $\square$ \bigskip}

\def\nkyear{202?}  
\def\nkvol{?}        
\def\nkno{?}          
\def\firstpage{1}      
\def\lastpage{pp}    
\def\correction{4}
\def\shorttitle{RBSDEs: monotone approximation} 
\def\shortauthor{{\it Bouhadou S, Hilbert A and Ouknine Y}} 
\def\shortauthor{{\it Bouhadou S,\ Hilbert A\ and\ \ Ouknine Y}} 
\usepackage{style}
\def\no#1#2#3{\parallel #1
	\parallel_{#2}^{#3}}

\def\f{{\cal F}}

\def\n{{\bf N}}

\def\ds{D_s}

\newcommand{\bdm}{\begin{displaymath}}
	\newcommand{\edm}{\end{displaymath}}
\newcommand{\bean}{\begin{eqnarray}}
	\newcommand{\eean}{\end{eqnarray}}
\newcommand{\bea}{\begin{eqnarray*}}
	\newcommand{\eea}{\end{eqnarray*}}

\def\BC{\mathcal{B}}
\def\FC{\mathcal{F}}
\def\PC{\mathcal{P}}

\def\R{{\bf R}}

\def\1B{\text{1\!\!I}}

\def\t{\tau}

\def\FC{\mathcal{F}}
\def\PC{\mathcal{P}}

\def\R{{\bf R}}

\def\1B{\text{1\!\!I}}

\def\HB{\mathbb{H}}

\newcommand{\cf}{\mathcal{F}}
\newcommand{\stopt}{\mathcal{T}_{t}}
\newcommand{\stops}{\mathcal{T}_{S}}
\newcommand{\stopo}{\mathcal{T}_{0}}

\newcommand{\barY}{\overline{Y}}

\title{\Large \bf \boldmath\ \\ Reflected backward stochastic differential equations with optional barriers: monotone approximation} 

\author{\large  Siham Bouhadou $^{1\dagger}$  Astrid Hilbert$^{2}$ and  Youssef Ouknine $^{ 1,3}$} 

\date{}

\begin{document}
	
	\maketitle
	
	\thispagestyle{first}
	\renewcommand{\thefootnote}{\fnsymbol{footnote}}
	
	\footnotetext{\hspace*{-5mm} \begin{tabular}{@{}r@{}p{13.4cm}@{}}
			& Manuscript received  \\ 
			$^1$ & Department of Mathematics, Faculty of Sciences Semlalia, Cadi Ayyad University,
			Morocco.\\
			&{E-mail:} sihambouhadou@gmail.com\\
			$^{2}$ &Department of Mathematics,  Linnaeus University, Sweden.\\
			$^{3}$ &Mohammed VI Polytechnic University, Morocco.\\
			$^{\dagger}$ & Corresponding author
	\end{tabular}}
	
	\renewcommand{\thefootnote}{\arabic{footnote}}
	
	\begin{abstract} 
		In this short note we consider RBSDE with Lipschitz drivers and barrier processes that are optional and right upper semicontinuous. We treat the case when the barrier can be represented as a decreasing limit of cadlag barriers. We combine well known existence results for cadlag barriers with comparison arguments for the control process to construct solutions. Finally, we highlight the connection of such RBSDEs with usual cadlag BSDEs.
		
		\vskip 4.5mm
		
		\nd \begin{tabular}{@{}l@{ }p{10.1cm}} {\bf Keywords } &
			reflected backward stochastic differential equation, g-expectation, optional barrier, monotone approximation,  comparison principle
		\end{tabular}
		
		\nd {\bf AMS Subject  Classification } 
		60H10; 60G40
	\end{abstract}
	
	\baselineskip 14pt
	
	\setlength{\parindent}{1.5em}
	
	\setcounter{section}{0}
	
	\Section{Introduction} \label{section1}
	
	Reflected backward stochastic differential equations (RBSDE) are a well known tool suited to solve the problem of hedging and pricing American options. The control process $Y$ of the solution triplet $(Y,Z,A)$ guiding the dynamics of RBSDE is reflected at a barrier process $\xi$, while the increasing process $A$ is responsible for keeping $Y$ above $\xi$. The original continuity assumption of El Karoui et al. \cite{Kar2} on $\xi$ has been relaxed in a series of papers to various degrees of discontinuity (see \cite{crepeymatoussi08, essaky08,hamadene02, hamadeneouknine03,hamadenewang09, hamadeneouknine11,klimsiak}). In \cite{Peng2}, Peng and Xu dealt with the case of a Brownain filtration and very irregular  $L^2$-obstacle, by introducing a new formulation of Skorokhod condition. The fundamental results on RBSDEs when the barrier is not right continuous were obtained in \cite{giooq15}. In \cite{BH} Bouhadou and Ouknine treated RBSDEs in the frame of a general filtration and a ladlag  predictable barrier. In these references cited above, the  construction of the solution was following the classical route of combining a priori inequalities with a recursively given sequence of approximations of the solution via a suitable fixed point argument.
	
	In this short note, we study the existence and uniqueness of the solution when the barrier can be approximated by a decreasing sequence of caldag Barriers $(\xi_n)_{n \in \mathbb{N}}$ in ${\cal S}^2$. This allows technically simpler proofs, since  we apply classical results on the existence of solution of RBSDEs with RCLL Barriers, and comparison related arguments (see Hamad\`ene, Wang \cite{hamadenewang09}) to show that the associated sequence of control processes $(Y^n)_{n\in\n}$ is also decreasing. Our first main result, shows that when the driver g does not depend on $y$ and $z$, the limiting process is the solution of the following optimal stopping problem.
	\begin{eqnarray}\label{eq20}	
		Y_S={\rm ess}\sup_{\tau\in \stops}E\left[\xi_\tau+ \int_S^\tau g(u) du \mid \FC_{S}\right].
	\end{eqnarray} 
	In order to prove that $Y$ provides the good candidate for the solution of RBSDEs with barrier $\xi$, we use some tools from optimal stopping theory (cf., \cite{Kar3}). \\ 
	In the second part of this note, we expand our convergence results to the non linear case, by considering the notion of $g$-conditional expectations (introduced by Peng \cite{Peng3}) defined through the notion of BSDEs  and used to quantify the riskiness of financial positions (see, among many others \cite{Delbaen,Barrieu,Peng3,Agnquenez}). We recall that the $g$-conditional expectation at a stopping time $\tau$ such that $\tau \leq T$ a.s.( where $T>0$ is a fixed final horizon) is the operator which maps a given square integrable terminal condition $\xi_T$ to the position at $\tau$ of the first component of the solution to the BSDEs with parameters $(g,\xi_T)$. The operator is denoted  by $\mathcal{E}^{g}(.)$.\\
	Roughly speaking, if we interpret $\xi$ as a financial position process and $-\mathcal{E}^{g}(.)$ as a dynamic risk measure, $\mathcal{R}^{g}[\xi](S)$
	defined by the following
	\begin{equation}
		\mathcal{R}^{g}[\xi](S):={\rm ess}\sup_{\t\in \mathcal{T}_S}\mathcal{E}^{g}_{S,\t}(\xi_\t),\;\;\;S \in \stopo,
	\end{equation}
	can be seen as the minimal risk at time $S$. In the present paper, we show that when the barrier can be approximated by RCLL barriers monotonically from above in $\mathcal{S}^2$,  the first component $Y$ of the solution of RBSDE with barrier $\xi$ satisfies the following
	$$Y_S=\mathcal{R}^{g}[\xi](S),\;\;\;S \in \stopo .$$
	In the last part, 
	we show how optional RBSDEs with non regular barrier,  are closely connected with usual RCLL BSDEs,  by introducing a new way to approximate the solution $Y$. The main novelty is to prove that the solution of RBSDE with optional barrier $\xi$ can be obtained as the limit of the following sequence processes:
	$$\bar Y_t^n=\xi_t \vee \left(\bar \xi_T+\int_t^Tg(u,Y_u^n,Z_u^n)du+\int_t^Tn(Y_u^n-\bar \xi_{u}^+)^-du-\int_t^T Z_u^n d W_u -\int_{t}^T\int_E  l_u^n(e) \tilde N(du,de) \right).$$
	Where $\bar \xi$ is given by  $\bar \xi=\mathcal{R}^{g}(\xi+X)-X$ and $X$ denotes an optional process satisfying some suitable assumptions.\\
	
	Let us present briefly our plan. In section $2$, we recall the solution  concept of optional RBSDE under a suitable version of Skorokhod condition. In section $3$, we present our solution by monotone approximation of the barrier. Section $4$, is dedicated to give a new approximation of the solution of RBSDE studied in section $2$. Section $5$, is devoted to make some link between RBSDE with optional Barrier and RBSDE in the sens of Peng-Xu \cite{Peng2}, under an additional assumption on $\xi$. 
	
	\section{Preliminaries} \label{sec2}
	First we introduce a series of notations that will be used throught the paper.
	Let $T>0$ be a fixed positive real number. Let $(E,\mathcal{E})$ be a measurable space equipped with a $\sigma$-finite positive measure $\mu$. Let $(\Omega, \mathcal{F}, P, \mathbb{F}=(
	\mathcal{F}_t)_{t\geq 0})$ be a probability space. The filtration is assumed to be
	complete, right continuous and quasi-left continuous.
	We suppose that $(\Omega, \mathcal{F}, P, \mathbb{F}=(
	\mathcal{F}_t)_{t\geq 0})$ supports a Brownian motion $W$ and an independent Poisson random measure $N$ with intensity $dt \otimes\mu(de)$. We denote $\tilde N(dt,de)$ its compensated Poisson random measure.\\
	For $t\in [0,T]$, we denote by $\stopt$ (resp. $\stopt^p$) the set of stopping times (resp. predictable stopping times)  $\tau$ such that $P(t \leq\tau\leq T)=1.$ More generally, for a given stopping time $\nu\in \stopo$  (resp.  $\nu\in \stopo^p$), we denote by $\mathcal{T}_{\nu}$ (resp. $\mathcal{T}_{\nu}^p$) the set of stopping times  (resp. predictable stopping times) $ \tau$ such that $P(\nu \leq\tau\leq T)=1.$ 	
	We denote by  ${\cal P}$ be  the predictable $\sigma$-algebra
	on $ \Omega\times [0,T]$. 
	We use the following notation: 
\begin{itemize}
\item $L^2({\cal F}_T)$  is the set of random variables which are  $\FC 
		_T$-measurable and square-integrable.
\item  $L^2_\mu$ is the set of measurable functions $\ell:  E \rightarrow \R$ such that  $\|\ell\|_\mu^2:= \int_{ E}  |\ell(e) |^2 \mu(de) < + \infty.$
\item  $\HB^{2}$ is the set of processes $\phi$ which are  {\em predictable} such that $$\| \phi \|^2_{\HB^{2}} :=E\left[ \int_0 ^T \phi_t^2 \,dt   \right]< \infty.$$
		
\item  $\HB_{\mu}^{2}$ is the set of processes $\phi$ which are  {\em predictable}, that is, measurable\\
$\phi: (\Omega\times [0,T] \times  E,\; \PC \otimes {\cal E}) \rightarrow (\R\;,  \BC(\R)); \quad 
		(\omega,t,e) \mapsto \phi_t(\omega, e)	$
		such that $$\| \phi \|^2_{\HB_{\mu}^{2}} :=E\left[ \int_0 ^T \|\phi_t\|_{\mu}^2 \,dt   \right]< \infty.$$
\item   $\mathbb{D}^{2}$ is the vector space of $\mathcal{F}_t$-adapted RCLL
		processes $\phi=(\phi)_{t\in[0,T] }$ such that
		$$\|\phi\|^2_{\mathbb{D}^{2}} := E[\sup_{t\leq T} |\phi_t |^2] <  \infty.$$
\item   ${\cal S}^{2}$ is the vector space of real-valued  optional 
		processes $\phi$ such that
		$$\vertiii{\phi}^2_{{\cal S}^{2}} := E[{\rm ess}\sup_{\tau\in\stopo} |\phi_\tau |^2] <  \infty.$$
		\item ${\cal S}^{4}$ is the vector space of real-valued  optional 
		processes $\phi$ such that
		$$\vertiii{\phi}^4_{{\cal S}^{4}} := E[{\rm ess}\sup_{\tau\in\stopo} |\phi_\tau |^4] <  \infty.$$
\item   ${\cal S}_{pi}^{2}$ is the vector space of real-valued  predictable, increasing 
		processes $A$ such that $A_0=0$, $E(A_T^2)<\infty$.
\end {itemize}
We say that an $\mathbb{F}$-progressively measurable process $X$ is of class $(D)$, if the family $\{X_\t, \;\t \in \stopo\}$ is uniformly integrable.

For a process $\psi$, we write $\psi_-$ for the process of left limits $\psi_{t-} = \lim_{s\uparrow t} \psi_s$, for $t>0$, provided they all exist, and $\psi^+$ for the process of right limits $\psi_{t+} = \lim_{s\downarrow t} \psi_s$ for $t<T$ in case they all exist. \\
For a ladlag process $X$, we denote by $\Delta^+ X_t:=X_{t_+}-X_t$ the size of the right jump of $X$ at $t$, and by   $\Delta X_t:=X_t-X_{t-}$ the size of the left jump of $X$ at $t$.\\ If $A$ is an increasing process, then it can be represented in the form $A=A^r+A^g$, with $A^r=A^c+A^d$, where $A^c$,  $A^d$ and  $A^g$ are increasing processes,  $A^c$ is a continuous component with $A_0^c=0$, $A^g$ is continuous from the left with $A_0^g=0$ and $A^d$ is continuous from the right with $\Delta A_0^d=A_0$.\\
Let us recall  the key section theorem related to indistinguishability of optional processes or predictable processes.
\begin{theorem}
Let $X=(X_t)$ and $Y=(Y_t)$ be two optional (resp. predictable) processes. If for every bounded stopping time (resp. predictable time) $\t$, we have $X_\t\leq Y_\t$ a.s. (resp. $X_{\tau}=Y_{\t}$ a.s.), then $X\leq Y$ (resp. X and Y are indistinguishable).
\end{theorem}
		\begin{definition}[Driver, Lipschitz driver]
			A function $g$ is said to be a {\em driver} if 
			\begin{itemize}
				\item  
				(measurability)	$g: \Omega  \times [0,T]\times \R^2 \times L^2_\mu \rightarrow \R $\\
				$(\omega, t,y, z, k) \mapsto  g(\omega, t,y, z, k)  $
				is $ {\cal P} \otimes {\cal B}(\R^2)  \otimes {\cal B}(L^2_\mu)$measurable,  
				\item (integrability) $ g(\cdot,0,0,0) \in \HB^2$.
			\end{itemize} 
			A driver $g$ is called a {\em Lipschitz driver} if moreover there exists a constant $ K \geq 0$ such that $dP \otimes dt$-a.e\,, 
			for each $(y_1, z_1, k_1) \in \R^2 \times L^2_\mu$ , $(y_2, z_2, k_2) \in \R^2 \times L^2_\mu$, 
			$$|g(\omega, t, y_1, z_1, k_1) - g(\omega, t, y_2, z_2, k_2)| \leq 
			K (|y_1 - y_2| + |z_1 - z_2| +\|k_1 - k_2 \|_\mu).$$
		\end{definition}
		Let  $g$ be a Lipschitz driver, and $\xi $ in $L^2(\cf_{T})$. The BSDE associated with Lipschitz driver $g$, terminal time $T$, and terminal condition $\xi$,  is formulated as follows:
		$$X_t=\xi+\int_{t}^T g(s, X_s, Z_s,l_s)ds-\int_{t}^T  Z_s dW_s-\int_{t}^T\int_E  l_s(e) \tilde N(ds,de)\;\; \text{for all } t\in[0,T] \text{ a.s. }  $$
		We recall that the above BSDE admits a unique solution $(X,Z,l)$ in the space $\mathbb{D}^2 \times \HB^{2} \times 	\HB_\mu^{2}$ (cf. \cite{essaky08}).\\
		We also recall the definition of the conditional $g$- expectation.
		
		\begin{definition}
			We define for each $t\in[0,T]$, and $\xi \in L^2(\cf_{T}) $ 
			
			$$\mathcal{E}^{g}_{t,T}(\xi):=X_t.$$
			
			We call  the (non-linear) operator $\mathcal{E}^{g}_{t,T}(\cdot): L^2(\cf_{T})\rightarrow L^2(\cf_{t})$
			\emph{conditional $g$-expectation at time $t$}.  As usual, this notion can be extended to
			the case where the (deterministic) terminal time 
			$T$
			is replaced by a (more general) stopping time $\tau\in\stopo$,  $t$ is replaced by a stopping time $S$ such that $S\leq \tau$ a.s. and  the domain $L^2(\cf_{T})$ of the operator    is replaced  by $L^2(\cf_{\tau})$.		
		\end{definition}
		Let $T>0$ be a fixed terminal time. Let $g$ be  a Lipschitz driver. 
		Let $\xi= (\xi_t)_{t\in[0,T]}$ be an optional  process in ${\cal S}^2$. 
		We suppose moreover that the process $\xi$ is not necessarily left limited. A process $\xi$ satisfying the previous properties will be called a  \emph{barrier}, or an  \emph{obstacle}.
		\begin{definition}\label{definition1}
			A quadruple $(Y,Z,l,A)$ of $\;\mathbb{F}$-progressively measurable processes is a solution of the reflected BSDE with Lipschitz driver $g$  and barrier $\xi$ $(RBSDE(\xi,g)$ for short) if
			\begin{itemize}
				\item[(i)] $(Y,Z,l,A)\in\mathcal{S}^{2}\times \HB^2 \times\HB^2_\mu \times \mathcal{S}^{2}_{pi} $.
				\item [(ii)] $Y_{\t}=\xi_T+\int_{\t}^{T}g(s,Y_{s},Z_s,l_s)ds+A_{T}-A_{\t}-\int_{\t}^T Z_s dW_s-\int_{\t}^T\int_E  l_s(e) \tilde N(ds,de)\!$ a.s. for all $\t\in\mathcal{T}_0$,
				\item[(iii)]  $Y_{\t} \geq \xi_{\t}$  a.s. for all $\t\in\mathcal{T}_0$,
				\item[(iv)]A is non decreasing predictable process with $A_0=0$ such that		
				$$\int_{0}^{T}(Y_{s-}-\limsup_{u \uparrow s} \xi_{u})dA_{s}^{r}=\sum_{s<T}(Y_{s}-\xi_{s})\Delta^{+}A_{s}=0 \;\;\;\text{a.s.}$$
				Here $A^r$ denotes the cadlag part of the process $A$.
			\end{itemize}
		\end{definition}
		
		\begin{remark}		
			Since the filtration is quasi-left continuous, martingales have only totally inaccessible jumps. Thus, in this case, ${}^{p}Y_\tau=Y_{\tau}$ for each predictable stopping time $\tau \in \stopo^p$.
		\end{remark}
		\begin{remark}\label{jumps-left}	
			It follows from (ii), and the same arguments as above that  $\Delta Y_{\tau}=-\Delta A^r_{\t}=-\Delta A_{\tau}$ a.s., for each predictable stopping time $\t$. This clearly yields that $Y_{\tau^-}\geq Y_{\tau} $ a.s. for each predictable stopping time $\t$.
		\end{remark}
		\begin{definition} \label{defr} A progressive process $(\xi_t)$ (resp. integrable) is said to be right (resp. left) upper semicontinuous   along stopping times (right (left) USC) (resp. along stopping times in expectation (right (left) USCE)) if for all $\tau \in {\cal T}_0$ and for all sequences of stopping times $ (\tau_n)$ such that  $\tau^n \downarrow \tau$ ( resp.  $\tau^n \uparrow \tau$)\,,
			\begin{eqnarray}\label{usc}
				\xi_{\tau} \geq \limsup_{n\to \infty} \xi_{\tau_n} \quad \mbox{a.s.} \quad \text{(resp.} E[\xi_{\tau}] \geq \limsup_{n\to \infty} E[\xi_{\tau_n}] ).
			\end{eqnarray}
		\end{definition}
		\begin{remark}\label{jumps}
			If $(Y,Z,l,A)$ is a solution of RBSDE defined above, then $\Delta^+ Y_{\tau}=Y_{\tau^+}-Y_{\tau}=-\Delta^+ A_{\tau}$ a.s. for each stopping time $\tau \in \stopo$.
			Roughly speaking, this equality says that the process has only negative right jumps.	
			Note also that $Y\geq Y^+$ up to an evanescent set, which means that process $Y$ is right upper semicontinuous.
		\end{remark}
		\begin{proposition}\label{prop-fon}
			Let $g$ be a Lipschitz driver and $\xi$ an obstacle. Let $(Y,Z,l,A)$ be a solution to the $RBSDE(\xi,g)$.
			\begin{itemize}
				\item  For each $\t \in \stopo$, we have 
				$$Y_\tau=\xi_\tau \vee Y_\tau^+\quad\mbox{a.s.} $$
				\item  For each predictable stopping time $\t \in \stopo^p$, we have 
				$$ Y_{\t^-}=\limsup_{u\uparrow \t}\xi_u \vee Y_{\t}\quad\mbox{a.s.} $$
			\end{itemize} 
		\end{proposition}
		\dproof
		Let us show the first assertion. Let  $\t \in \stopo$. The inequality
		$\xi_\tau \vee Y_\tau^+ \leq Y_\tau$ a.s. follows from the fact that $\xi_\tau \leq Y_\tau $ a.s. and $Y_\tau^+ \leq Y_\tau$  a.s. Let us now show the second inequality.   Thanks to Remark \ref{jumps}, $\Delta^+A_\t=-\Delta^+ Y_\t$ a.s. Then, from Skorokhod condition $(iv)$, we have  $Y_\tau 1_{\{Y_\tau>\xi_\tau\}}= Y_\tau^+ 1_{\{Y_\tau>\xi_\tau\}}$ a.s. that $Y_\tau\leq \xi_\tau \vee Y_\tau^+$ a.s. \\The task now is to prove the second assertion. Let $\t \in \stopo^p$. We have  $\limsup_{u\uparrow \t}\xi_u \leq Y_{\t^-}$ a.s. and $Y_{\t} \leq Y_{\t^-}$ a.s. Hence, $\limsup_{u\uparrow \t}\xi_u \vee Y_{\t} \leq Y_{\t^-}$ a.s. Let us now focus on the first inequality. It follows from Remark \ref{jumps-left}, that $\Delta Y_{\tau}=-\Delta A^r_{\t}$. Then, using the second inequality of the Skorokhod condition we get,  $Y_\t 1_{\{Y_{\t^-}>\limsup_{u\uparrow \t}\xi_u\}}=Y_{\t^-} 1_{\{Y_{\t^-}>\limsup_{u\uparrow \t}\xi_u\}}$ a.s. Thus we have proved that $ Y_{\t^-} \leq \limsup_{u\uparrow \t}\xi_u \vee Y_{\t}$ a.s.
		\fproof
		\begin{lemma}\label{lem left}If $\xi$ is left USC, then $Y_{\tau}=Y_{\tau^-}$ for each predictable stopping time $\tau \in \stopo^p$. On other words, the process $A^r$ is continuous.
		\end{lemma}
		\dproof
		Let $\tau$ be a predictable stopping time in  $\stopo^p$.
		The second assertion in Proposition \ref{prop-fon}, combined with the fact that $\xi$ is left USC leads to $Y_{\tau^-}\leq \xi_\tau\vee Y_{\tau}  =Y_{\tau}$ a.s. Otherwise we know that $Y_{\tau}  \leq Y_{\tau^-}$ a.s,  at last $Y_{\tau}=Y_{\tau^-}$ a.s. The continuity of the process $A^r$ follows immediately from Remark \ref{jumps-left}.
		\fproof
\begin{remark}
If $Y$ is right continuous, then the Skorokhod condition (iv), can be reduced to the following:
$$\int_{0}^{T}(Y_{s-}-\limsup_{u \uparrow s} \xi_{u})dA_{s}^{r}=0 \;\;\;\mbox{a.s.}$$
Indeed, the right continuity of $Y$ together with the Remark \ref{jumps} implies that $\Delta^+A_t=0$ a.s. for all $t \in [0,T]$.	
\end{remark}
		\begin{lemma}\label{lem right}
			If the obstacle $\xi$ satisfies $\xi\leq \xi^+$ up to an evanescent set, then $Y$ is right-continuous.
		\end{lemma}
		\dproof
		Through the first assertion of Proposition \ref{prop-fon}, we have for any $\tau \in \stopo$, $Y_\tau=\xi_\tau \vee Y_\tau^+ $ a.s. Under the assumption on $\xi$, we obtain $Y_\tau=\xi_\tau \vee Y_{\tau^+} \leq \xi_{\tau^+} \vee Y_{\tau^+} =Y_{\tau^+}$ a.s. Thanks to  Proposition \ref{jumps},
		$Y_{\tau^+} \leq Y_{\tau}$ a.s. This ends the proof.
		\fproof	
		
		Next we introduce the notion of strong supermartingale which extend the classical supermartingales to those connected to  the optional $\sigma-$field.
		\begin{definition}An optional process $(Y)_{t \in [0,T]}$ such that
			\begin{itemize} 
				\item $Y_\tau $ is integrable for all $\tau \in \mathcal{T}_0$.
				\item for arbitary stopping times  $\tau \geq \sigma$
				$$Y_\sigma\geq E[Y_\tau\mid \FC_{\sigma}]\;\;\;\;\mbox{a.e}$$ 
				is called a strong supermartingale.
			\end{itemize}
		\end{definition}
		\begin{remark}\label{rem.ladlag}
			Every optional strong supermartingale is indistinguishable from a ladlag process,  see \cite{DM2}.
		\end{remark}
		\section{Monotone approximation of the barrier}\label{s:monotone}
		\begin{definition}\label{def:lower-envelope}
			Let $\xi$ be an optional process. Let
			$$\mathcal{L} = \{ X: X \,\,\mbox{is a cadlag optional process,}\,\,X\ge \xi\},\quad \mathcal{L}_- = \{X_-: X \in\mathcal{L}\},$$
			and
			$$\overline{\xi} = \mbox{ess inf}\,\, \mathcal{L},\quad\quad\quad\quad  \hat{\xi} = \mbox{ess inf}\,\, \mathcal{L}_-.$$
			We call $\overline{\xi}$ \emph{upper cadlag envelope} of $\xi$, $\hat{\xi}$ \emph{left upper cadlag envelope} of $\xi$.
		\end{definition}
		Let us provide some properties of the process $\bar \xi$.
		\begin{lemma}
			Let $\xi$ be an optional process. Then, there exists a non increasing sequence $(X^n)_{n\in\n}$ of cadlag processes in $\mathcal{L}$ such that
			\begin{eqnarray}
				\overline{\xi} &=& \lim_{n\to\infty}\downarrow X^n.
			\end{eqnarray}		
		\end{lemma}
		\dproof
		The  infima of a finite number of processes in $\mathcal{L}$ belongs to $\mathcal{L}$. Thus by a result of Neveu \cite{neveu}, $\mbox{ess inf}\,\,\mathcal{L}$ can be described as the infimum of a sequence of processes in $\mathcal{L}$. This concludes the proof.
		\fproof
		\begin{remark}
			It is clear that if the optional process $\xi$ is RCLL then, its upper cadlag envelope is more closely related to $\xi$. We give a generalization of this result to the case of right upper semicontinuous process in the following lemma:
		\end{remark}
		
		\begin{lemma}\label{lem:envelope}
			Let $\xi$ be an optional process which is right upper semicontinuous. Then
			$$\overline{\xi}_t = \xi_t, \quad t\in[0,T],$$
			i.e. $\overline{\xi}$ is a version of $\xi.$ Moreover, the sequence $(X^n)_{n\in\n}$ may be chosen identical to the sequence $(\xi^n)_{n\in\n}$ resulting from Theorem 21 of Dellacherie, Lenglart \cite{dellacherielenglart82}. Finally, we have
			$$\hat{\xi} = \lim_{n\to\infty} {\xi^n_-}.$$
		\end{lemma}
		
		{\bf Proof:}\\
		Since for any $n\in\n$,  $\xi^n$ is RCLL and optional, we have for $t\in[0,T]$:
		$$\overline{\xi}_t \le \lim_{n\to\infty} \xi^n_t = \xi_t\le \overline{\xi}_t.$$
		Thus the equation on $\hat{\xi}$ follows from the definition of $\hat{\xi}.$
		\hfill $\Box$
		\begin{remark}
			If  $\xi$ is an optional process which is right upper semicontinuous in expectation of class $(D)$, then it is right upper semicontinuous.	
		\end{remark}	
		\begin{lemma}\label{expectation1}
			Let $\xi$ be an optional process which is right upper semicontinuous in expectation. Then
			$$\overline{\xi}_t = \xi_t, \quad t\in[0,T],$$
		\end{lemma}
		In particular, we get the following lemma:
		\begin{lemma}\label{supermartingale}
			Let $S$ be an optional strong supermartingale of class $(D)$. Then
			$$\overline{S}_t = S_t, \quad t\in[0,T],$$
		\end{lemma}
		\dproof
		Since $S$ is an optional strong supermartingale, the application $\tau \rightarrow E(S_\tau)$ is non increasing. Thus, $S$ is clearly right upper semicontinuous in expectation. The result follows from an application of Lemma \ref{expectation1}.
		\fproof
		\subsection{Optional RBSDEs from RCLL RBSDEs }
		We now show how the solution of the RBSDE when optional barrier is approximated by RCLL barriers monotonically from above in $\mathcal{S}^2$,  can alternatively be constructed along a sequence of RBSDE with RCLL barriers. Let $g$ be a Lipschitz driver.\\
		Assume that $\xi$ is right upper semicontinuous. Through  the proof of  Proposition $21$ in Dellecherie-Lenglart \cite{dellacherielenglart82}, there exits a sequence $(\xi^n )_{n \in \mathbb{N}}\in\mathcal{L}$ such that  $E[\sup_{t \in [0,T]}(\xi_t^n)^2] < \infty$ and 
\begin{eqnarray}\label{approx}
\xi^n \downarrow  essinf \mathcal{L}=\xi.
\end{eqnarray}
		Moreover, we assume that $\|\xi^n-\xi\|_{\mathcal{S}^2}\rightarrow 0$ as $n \rightarrow \infty$.
		By Hamad\`ene and  Ouknine \cite{hamadeneouknine11}, there exists  $(Y^n, Z^n,l^n,A^n)$  the solution of the following RCLL reflected BSDE:
		\begin{equation}
			\left\{
			\begin{array}{ll}
				(i)& (Y^n, Z^n,l^n,A^n)\in \mathbb{D}^2 \times \HB^2\times  \HB^2_{\mu}\times {\cal S}_{pi}^2,
				\nonumber\\
				(ii)& Y_t^n=\xi_T^n + \int_t^T g(s, Y_s^n, Z_s^n,l_s^n) ds - \int_t^T Z_s^n d W_s -\int_{t}^T\int_E  l_s^n(e) \tilde N(ds,de) + A ^n_T - A_t^n,\;\mbox{for all}\; t\in[0,T], \;\mbox{a.s.} \\\label{prince} 
				(iii)&  Y_t^n \geq \xi^n_t, \;\mbox{for all} \;t \in[0,T] \;\mbox{a.s.},\\ \nonumber
				(iv)& A^n\,\, \mbox{is cadlag predictable, increasing with}\,\, A^n_0=0, E(A^n_T)<\infty, \,\,\mbox{and satisfies}\\
				& \int_0^T (Y_{t-}^n - \xi_{t-}^n) dA_t^n = 0. 
			\end{array}
			\right.
		\end{equation}
		
		Note that $dA^n$ is the (random) measure on the Borel sets of $[0,T]$ associated with the increasing cadlag process $A^n$. We also remark that the {\em Skorokhod condition} in ($iv$) can be translated into the more detailed condition:
		
		\begin{eqnarray}\label{eq:Skorokhod-n-alternative}
			\int_0^T (Y_t^n-\xi_t^n) d A^{n,c}_t = 0,\quad \int_0^T (Y_{t-}^n - \xi_{t-}^n) dA_t^{n,d} = 0.
		\end{eqnarray}
		
		Here $A^{n,c}$ denotes the continuous part of $A^n$, $A^{n,d}$ its discontinuous part.
		\begin{theorem}\label{thm:existence-n}
			Assume that $\xi$ is in $\mathcal{S}^2$. Let a sequence of decreasing RCLL processes $(\xi^n)_{n\in\n}$ be given which satisfies (ii). Let $g$ be a Lipschitz driver. Then, for each $n\in\n$,  there exists a quadruple  $(Y^n, Z^n,l^n,A^n)\in \mathbb{D}^2\times \mathbb{H}^2 \times\mathbb{H}^2_{\mu}\times\mathcal{S}_{pi}^2$ of processes which solves the RBSDE (ii), (iii), (iv). Moreover, we have for any $n\geq 0$ and for any $t\in[0,T]$:
			\bean
			Y_t^n &\ge& Y_t^{n+1}.\label{eq:domination-solution-n}
			\eean
		\end{theorem}
		\dproof
		This is shown in Hamad\`ene, Wang \cite{hamadenewang09}.
		\fproof
		
		Let us now investigate the convergence of the first component of the solution quadruple of Theorem \ref{thm:existence-n}. First of all, for $n\in\n, t\in[0,T],$ we have
		$$\xi_t\leq Y_t^n \leq Y_t^1.$$
		Hence, by square integrability of $Y^1$ and the fact that $\xi\in\mathcal{S}^2$, and by dominated convergence, the sequence $(Y^n)_{n\in\n}$ converges in $L^2(\Omega\times[0,T], \mathcal{F}_T\otimes \mathcal{B}([0,T]), P\otimes dt)$ to a process $Y\in\mathcal{S}^2$. 
		The proof of our main result is based on the following key theorem:
		\begin{theorem}\label{thm:main1}
			Suppose that $g$ does not depend on $y$, $z$, $l$ that is $g(\omega,t,y,z,l)=g(\omega,t)$,   where $g$ is a progressive process with $E(\int_0^t g(t)^2 )dt<+ \infty$.
			Let $\xi$ be an optional process which is right upper semicontinuous, and let $(\xi^n)_{n\in\n}$ be given according to (\ref{approx}). Let $Y = \lim_{n\to\infty} Y^n$ in $L^2(\Omega\times[0,T], \mathcal{F}_T\otimes \mathcal{B}([0,T]), P\otimes dt)$,
			then, for each $S \in \stopo$
			\begin{equation}\label{eq1}	
				Y_S= {\rm ess}\sup_{\tau\in \stops}E\left[ \xi_\tau+ \int_S^\tau g(u) du \mid \FC_{S}\right].
			\end{equation} 
			And the following  properties hold:
			\begin{description}
				\item[(i)]	We have	$ Y\equiv\xi\vee Y^{+}$ .
				\item[(ii)]	We have	$ Y_{\t^-}=\limsup_{u\uparrow \t}\xi_u \vee Y_{\t}$, for all $\t\in\stopo^p$.
			\end{description}
			Moreover,  the convergence of the sequence $(Y^n)_{n\in\n}$ to $Y$  holds in $\mathcal{S}^2$.
		\end{theorem}
		\dproof
		First, let us show the equality (\ref{eq1}). We get from  Hamad\`ene and Ouknine \cite{hamadeneouknine11}, that for $ n\geq 0$ and $S \in \stopo$:
		\begin{eqnarray*}	
			Y_S^n= {\rm ess}\sup_{\tau\in \stops}E\left[ \xi_\tau^n+ \int_S^\tau g(u) du \mid \FC_{S}\right].
		\end{eqnarray*} 
		Let $\sigma \in \stopo$, 	Let  us denote  $\bar Y(\sigma)$ the random variable defined by: $$\bar Y(\sigma):={\rm ess}\sup_{\tau\in \mathcal{T}_\sigma}E\left[ \xi_\tau+ \int_\sigma^\tau g(u) du \mid \FC_{\sigma}\right].$$
		Therefore,
		$$\bar Y(\sigma)+\int_0^\sigma g(u)du ={\rm ess}\sup_{\tau\in \mathcal{T}_\sigma}E\left[ \xi_\tau+ \int_0^\tau g(u) du \mid \FC_{\sigma}\right].$$
		Since the process $(\xi_\cdot+ \int_0^\cdot g(u) du)$ is of class $(D)$, the family $(\bar Y(\sigma),\;\sigma \in \stopo)$ can be aggregated by a process wich we denote also  $\bar Y$ (cf., \cite[Theorem 15 ]{dellacherielenglart82}). Therefore
		
		\begin{eqnarray}\label{pro}
			&|\barY_\sigma-Y_\sigma^n|\leq {\rm ess}\sup_{\t\in \mathcal{T}_\sigma}E\left[| \xi_\tau^n-\xi_\tau| \mid \!\!\FC_{\sigma}\right]\\\nonumber
			&\leq E\left[ {\rm ess}\sup_{\t\in \mathcal{T}_\sigma}| \xi_\tau^n-\xi_\tau|\mid \FC_{\sigma}\right].\nonumber
		\end{eqnarray}
		First, let us denote $$U_\sigma^n=E\left[{\rm ess}\sup_{\tau \in \mathcal{T}_0}|\xi_\tau^n-\xi_\tau| \mid \FC_{\sigma}\right].$$ 
		
		Note that the process $(U_t^n)_{t \in [0,T]}$ is right continuous. This together with the definition of the essential supremum give $${\rm ess}\sup_{\sigma \in \mathcal{T}_0}|U_\sigma^n|^2=\sup_{t \in [0,T] } |U_t^n|^2\;\;\;\; \mbox{a.s.}$$
		By using (\ref{pro}), we obtain
		\begin{equation*}
			{\rm ess}\sup_{\sigma\in \mathcal{T}_0}	|\barY_\sigma-Y_\sigma^n|^2\leq\sup_{t \in [0,T] } |U_t^n|^2\;\;\; \mbox{a.s.}\nonumber
		\end{equation*}
		We apply  Doob's inequality, to get
		$$\|\barY-Y^n\|^2_{\mathcal{S}^{2}}\leq 
		E\left[{\rm ess}\sup_{\tau \in \mathcal{T}_0}|\xi_\tau^n-\xi_\tau|^2\right].$$
		The sequence  $(\xi^n)_{n \in \mathbb{N}}$ converges to $\xi$ in $\mathcal{S}^2$ by hypothesis. 
		Therefore, $\|\barY-Y^n\|_{\mathcal{S}^2}\rightarrow 0$ as $n \rightarrow \infty$. 
		Let us recall that  the sequence $(Y^n)_{n\in\n}$ converges in $L^2(\Omega\times[0,T], \mathcal{F}_T\otimes \mathcal{B}([0,T]), P\otimes dt)$ to a process $Y\in\mathcal{S}^2$. 
		Thus, we get $$\bar Y(\sigma)=\bar Y_\sigma= Y_\sigma\quad \mbox{ a.s. for all} \;\;
		\sigma \in \stopo.$$
		Which establishes that the process  $(Y_t+\int_0^tg(u)du)_{t \in [0,T]}$ is indistinguishable from the Snell envelope of the process $
		(\xi_t+\int_0^tg(u)du)_{t \in [0,T]}$. Assertions (i) and (ii) follow from classical  results (cf., for instance \cite[Proposition 2.32]{Kar3}).\\
		\fproof
		\begin{theorem}\label{princ}
			Let $\xi$ be a right upper semicontinuous barrier which can be approximated by RCLL barriers monotonically from above in $\mathcal{S}^2$. Suppose that $g$ does not depend on $y$, $z$, $l$ that is $g(\omega,t,y,z,l)=g(\omega,t)$,   where $g$ is a progressive process with $E(\int_0^t g(t)^2 )dt<+ \infty$. The reflected BSDE with one reflecting barrier associated with $(g,\xi)$  has a unique solution $(Y,Z,l,A)$. Where  $Y$ is given according to  Theorem \ref{thm:main}. Moreover, 
			the first component can be characterized as follows:
			\begin{equation}	
				Y_S= {\rm ess}\sup_{\tau\in \stops}E\left[ \xi_\tau+ \int_S^\tau g(u) du \mid \FC_{S}\right],\;\mbox{for all}\;S\in \stopo,
			\end{equation}
			and the following  properties hold:
			\begin{description}
				\item[(i)]	We have	$ Y\equiv\xi\vee Y^{+}$ .
				\item[(ii)]	We have	$ Y_{\t^-}=\limsup_{u\uparrow \t}\xi_u \vee Y_{\t}$ \;\mbox{a.s.},\; for all $\t\in\stopo^p$.
			\end{description}
		\end{theorem}
		The proof of Theorem \ref{princ} relies on the Theorem \ref{thm:main1} and the following lemma:
		\begin{lemma}\label{Mertensvalue}
			\begin{description}
				\item[(i)]  The  process $( Y_t)_{t\in[0,T]}$ is in ${\cal S}^{2}$ and admits the following optional Mertens decomposition:
				\begin{equation}\label{eqmert}
					Y_{\tau}=Y_{0}-\int_0^\tau g(s,Y_s,Z_s)ds +\int_0^T Z_s d W_s +\int_{0}^T\int_E  l_s(e) \tilde N(ds,de)-A_{\tau},\; \mbox{for all}\;\tau\in\stopo.
				\end{equation}
				where $A$ is a nondecreasing
				optional process such that $ A_0= 0$ and	$E(A_T^2)<\infty$.
				\item[(ii)] 
				$\int_{0}^{T}(Y_{s-}-\limsup_{u\uparrow s}\xi_u)dA_{s}=\sum_{s<T}(Y_{s}-\xi_{s})\Delta^{+}A_{s}=0$ a.s.
			\end{description}	
		\end{lemma}
		\dproof
		Let us prove the first assertion. For each $S \in \stopo$, we  define the random variable $U(S)$ by 
		$$U(S):=Y_S+\int_0^S g(u)du={\rm ess}\sup_{\t \in \stops} E\left[ \xi_\tau + \int_0^\t g(u)du \mid \mathcal{F}_{S}\right].$$
		By \cite{kob}, the process $(Y_t+\int_0^tg(u)du)_{t \in [0,T]}$ is the Snell's envelope associated to $(\xi_t + \int_0^t g(u)du)$. By this and by using Merten's decomposition, we get the equation (\ref{eqmert}).\\
		Now, let us show the assertion $(ii)$.\\First let us note that:
		$$\int_{0}^{T}(Y_{s-}-\limsup_{u\uparrow s}\xi_u)dA_{s}=0,$$
		can be written as the following:
		$$\int_{0}^{T}(Y_{s-}-\limsup_{u\uparrow s}\xi_u)dA_{s}^c=0\;\mbox{a.s.}\;\;\;\; \sum_{s \leq T}(Y_{s-}-\limsup_{u\uparrow s}\xi_u)\Delta A_{s}=0\;\mbox{a.s.}$$
		
		The proof  of the first inequality is based on the same arguments used in \cite{grig2}.\\
		The second equality is a consequence of (ii) in Theorem \ref{thm:main1} and Remark \ref{defr}. The following equality
		$$\sum_{s<T}(Y_{s}-\xi_{s})\Delta^{+}A_{s}=0,$$
		follows from (i) in Theorem \ref{princ} 
		and Remark \ref{jumps}.
		\fproof
		\section{Existence and uniquenes in the case of a general driver}
		Let  $Y^n$ be the first component of the solution of the RBSDE with the  RCLL barrier $\xi^n$ and the driver $g$. In \cite{quen1}, the authors proved that the vlaue function of the optimal stopping problem coincides with $Y^n$. Roughly speaking, for each stopping time $S \in \mathcal{T}_0$,
		
		$$Y_S^n={\rm ess}\sup_{\tau\in\mathcal{T}_S}\mathcal{E}^{g}_{S,\tau}(\xi_\t^n).$$
		
		In this part, we will use this characterisation, to construct the solution of the RBSDE when the optional barrier $\xi$ can be approximated by RCLL barriers $ \xi^n$ monotonically from above in $\mathcal{S}^2$.\\
		But, first les us revisit some properties of the family of random variables $(\mathcal{R}^{g}[\xi](S), \; S\in\stopo)$ defined by:
		
		\begin{equation}\label{nonlinear}
			\mathcal{R}^{g}[\xi](S):={\rm ess}\sup_{\t\in \mathcal{T}_S}\mathcal{E}^{g}_{S,\t}(\xi_\t),\;\;\;S \in \stopo.
		\end{equation}
		Let us recall the following definition:
		\begin{definition}\label{def.admi}
			We say that a family 
			$\phi=(\phi(\theta), \, \theta\in \mathcal{T}_0)$  is \emph{admissible} if it satisfies the following conditions 
			\par
			1. \quad for all
			$\theta\in \mathcal{T}_0$, $\phi(\theta)$ is a $\mathcal{F}_\theta$-measurable random variable,
			\par
			2. \quad  for all
			$\theta,\theta'\in \mathcal{T}_0$, $\phi(\theta)=\phi(\theta')$ a.s.  on
			$\{\theta=\theta'\}$.\\
		\end{definition}
		\begin{definition}An admissible square-integrable family $U:=(U(\theta), \; \theta \in\stopo)$ is
			said to be a strong  ${\cal E}^{g}$-\emph{supermartingale family} (resp. a strong ${\cal E}^{g}$-\emph{martingale family}), if for any 
			$\theta, \theta^{\prime}$ $ \in$ $\mathcal{T}_0$ such that $\theta^{\prime} \geq \theta$ a.s.,
			\begin{equation}
				{\cal E}^{g}_{\theta,\theta^\prime} (U(\theta^\prime))  \leq  U(\theta) \quad \,\mbox{a.s.}
				{\rm (resp.} \quad 	{\cal E}^{g}_{\theta,\theta^\prime} (U(\theta^\prime)) =  U(\theta)).
			\end{equation}
		\end{definition} 
		The following proposition plays a crucial role to derive some properties of the family $(\mathcal{R}^{g}[\xi](S), \; S\in\stopo)$.
		\begin{proposition}\label{prop.bellman}
			Let $S \leq  \theta \in \stopo$, and  let $\alpha $ be a non negative bounded $\mathcal{F}_{\theta}$-measurable random variable. Then,
			\begin{eqnarray}\label{bellman}
				\mathcal{E}^{\alpha g}_{S,\theta}(\alpha \mathcal{R}^{g}[\xi](\theta)) ={\rm ess}\sup_{\theta \leq \tau \in \stopo}\mathcal{E}^{\alpha g}_{S,\tau}(\alpha\xi_\t)\;\mbox{a.s.}
			\end{eqnarray}
		\end{proposition}
		
		\dproof
		Let $\tau \in \mathcal{T}_\theta$. By using the consistency property of $g$-conditional expectations, the fact that $\alpha$ is $\mathcal{F}_\theta$-measurable, $\mathcal{E}^{g}_{\theta,\tau}(\xi_\tau) \leq  \mathcal{R}^{g}[\xi](\theta)$ a.s. and the monotonicity  property of $g$-conditional expectations, we obtain:
		$$\mathcal{E}^{\alpha g}_{S,\tau}(\alpha\xi_\tau)=\mathcal{E}^{\alpha g}
		_{S,\theta}(\mathcal{E}^{\alpha g}_{\theta,\tau}(\alpha \xi_\tau))=\mathcal{E}^{\alpha g}_{S,\theta}(\alpha\mathcal{E}^{g}_{\theta,\tau}( \xi_\tau))\leq \mathcal{E}^{\alpha g}_{S,\theta}(\alpha \mathcal{R}^{g}[\xi](\theta))$$
		By taking the essential supremum over
		$\tau \in \mathcal{T}_\theta$, the  inequality:
		$$ {\rm ess}\sup_{ \theta \leq \tau \in \stopo}\mathcal{E}^{\alpha g}_{S,\tau}(\alpha \xi_\tau)\leq \mathcal{E}^{\alpha g}_{S,\theta}(\alpha \mathcal{R}^{g}[\xi](\theta))\quad
		\mbox{\rm a.s.}$$
		holds.
		We need to show the reverse inequality. Following \cite{giooq15}, there exists a sequence $(\tau_n)_{n\in \mathbb{N}}$  of stopping times in $\mathcal{T}_\theta$ such that the
		sequence $({\cal E}^{g}_{\theta,\tau_n}(\xi_{\tau_{n}}))_{n \in \mathbb{N}}$ is
		non decreasing and:
		$$\mathcal{R}^{g}[\xi](\theta)=\lim_{n \to \infty} \uparrow {\cal E}^{g}_{\theta,\tau_n}(\xi_{\tau_{n}})\quad
		\mbox{\rm a.s.}$$
		
		By using the fact that $\alpha$ is $\mathcal{F}_{\theta}$-measurable and a standard property of conditional $g$-expectations, (cf., e.g., Proposition $2.2$ in \cite{Peng3}), we obtain:  $$\alpha \mathcal{R}^{g}[\xi](\theta)=\lim_{n \to \infty} \uparrow {\cal E}^{\alpha g}_{\theta,\tau_n}(\alpha\xi_{\tau_{n}})\;
		\mbox{\rm a.s.}$$
		Therefore, by applying the property of continuity of BSDEs with respect to terminal condition (cf., Proposition A.$6$  in \cite{Agnquenez}) combined with consistency property of $g$-conditional expectations, we get
		$$\mathcal{E}^{\alpha g}_{S,\theta}(\alpha \mathcal{R}^{g}[\xi](\theta)) =\lim_{n \to \infty} \mathcal{E}^{\alpha g}_{S,\theta}({\cal E}^{\alpha g}_{\theta,\tau_n}(\alpha\xi_{\tau_{n}}))=\lim_{n \to \infty} {\cal E}^{\alpha g}_{S,\tau_n}(\alpha\xi_{\tau_{n}}).$$
		Hence,
		$$\mathcal{E}^{\alpha g}_{S,\theta}(\alpha \mathcal{R}^{g}[\xi]( \theta)) \leq {\rm ess}\sup_{\theta \leq \tau \in \stopo}\mathcal{E}^{\alpha g}_{S,\tau}(\alpha\xi_\tau).$$
		Whence the desired result.
		\fproof
		\begin{proposition}\label{prop.SuperM} 
			The value  family  $(\mathcal{R}^{g}[\xi](S), S\in \mathcal{T}_0)$ is characterized as the \emph{strong ${\cal E}^{g}$- Snell envelope family}  associated with $\xi$, that is,  the smallest  ${\cal
				E}^{g}$-supermartingale family which is greater (a.s.) than  or equal to $\xi$.
		\end{proposition}
		\dproof Let $\theta \in \mathcal{T}_S$.
		Applying the Proposition \ref{prop.bellman} with $\alpha=1$, and using that $S \leq  \theta$ a.s. we get:
		$$\mathcal{E}^g_{S,\theta}(\mathcal{R}^{g}[\xi]( \theta)) ={\rm ess}\sup_{ \tau \in \mathcal{T}_{\theta}}\mathcal{E}^g_{S,\tau}(\xi_\tau)\leq{\rm ess}\sup_{ S \leq \tau \in \stopo}\mathcal{E}^g_{S,\tau}(\xi_\tau) =\mathcal{R}^{g}[\xi](S) \quad
		\mbox{\rm a.s.}$$
		It follows that $\mathcal{R}^{g}[\xi]$ is an ${\cal E}^{g}$-supermartingale family. To complete the proof, it remains to show the minimality property. Let $V^\prime$ another	${\cal
			E}^{g}$-supermartingale family, such that $V^\prime\geq \xi$. The monotonicity property of $g$-conditional expectations allows us to write:
		$$\mathcal{E}^g_{S,\theta}(\xi_\theta)\leq \mathcal{E}^g_{S,\theta}(V^\prime(\theta))\leq V^\prime(S) \quad
		\mbox{\rm a.s.}, $$
		where the last inequality is due to the $\mathcal{E}^g$-supermartingale property of $V'$. By taking the essential supremum over
		$\theta \in \mathcal{T}_S$, we deduce that
		
		$$\mathcal{R}^{g}[\xi](S)={\rm ess}\sup_{\theta\in \mathcal{T}_S}\mathcal{E}^{g}_{S,\theta}(\xi_\theta) \leq V'(S) \quad
		\mbox{\rm a.s.}$$
		This concludes the proof.
		\fproof
		\begin{remark}\label{remark3}
			If $\xi \leq \tilde \xi$, then, $\mathcal{R}^{g}[\xi](S) \leq \mathcal{R}^{g}[\tilde\xi](S)$ a.s. for all $S \in \stopo$. First let us notice that through the defintion of $g$-conditional expectations and comparison theorem for BSDEs (cf. for e.g. \cite{essaky08}), we get for all $\t \in \mathcal{T}_{S}$
			$$\mathcal{E}^g_{S,\t}(\xi_\t)\leq \mathcal{E}^g_{S,\t}(\tilde \xi_\tau)\;\;\mbox{a.s.}$$
			We conclude the inequality by taking the essential supremum over $\t \in \mathcal{T}_S$.
		\end{remark}
		The purpose of the following theorem, is to  show that under suitable type of convergence of a sequence of reward processes $(\xi^n)_{n \in \mathbb{N}}$, the following convergence in terms of BSDEs can be proved, by using some a priori estimates of BSDEs.
		
		\begin{theorem}\label{cornonliear}
			Let $\xi$ be an optional process in $\mathcal{S}^2$. Let $(\xi^n)_{n \in \n}$ be a sequence of  optional processes in $\mathcal{S}^2$ such that $\|\xi^n-\xi\|_{\mathcal{S}^2} \rightarrow 0$. Then, the sequence $ (\mathcal{R}^{g}[\xi^n])_{n \in \mathbb{N}}$  converges in $\mathcal{S}^4$ to  $\mathcal{R}^{g}[\xi]$.
		\end{theorem}
		\dproof
		We have: 
		\begin{eqnarray}\label{general-est}
			{\rm ess}\sup_{ \theta  \in \stopo}|\mathcal{R}^{g}[\xi^n](\theta)-\mathcal{R}^{g}[\xi](\theta)|^4&=&	{\rm ess}\sup_{\theta\in \stopo}|{\rm ess}\sup_{\tau \in \mathcal{T}_{\theta}}\mathcal{E}^{g}_{\theta,\tau}(\xi_{\tau}^n)-{\rm ess}\sup_{ \tau \in \mathcal{T}_{\theta}}\mathcal{E}^{g}_{\theta,\t}(\xi_\t)|^4\\\nonumber
			\nonumber &\leq &{\rm ess}\sup_{\theta\in \stopo}{\rm ess}\sup_{\tau \in \mathcal{T}_{\theta}}|\mathcal{E}^{g}_{\theta,\t}(\xi_{\t}^n)-\mathcal{E}^{g}_{\theta,\t}(\xi_{\tau})|^4.
		\end{eqnarray}
		On the other hand, we have by a priori  estimates of BSDEs  (cf., Proposition A.$6$  in \cite{Agnquenez}), for each $\tau \in \mathcal{T}_\theta$
		\begin{eqnarray}
			|\mathcal{E}^{g}_{\theta,\t}(\xi_{\t}^n)-\mathcal{E}^{g}_{\theta,\t}(\xi_{\tau})|^4&\leq& c(E[|\xi_{\t}^n-\xi_{\t}|^2|\mathcal{F}_{\theta}])^2.\nonumber
		\end{eqnarray}
		Here $c$ is a constant which can changes from line to line.\\
		Thus,
		\begin{eqnarray}
			{\rm ess}\sup_{\theta\in \stopo}{\rm ess}\sup_{\tau \in \mathcal{T}_{\theta}}|\mathcal{E}^{g}_{\theta,\t}(\xi_{\t}^n)-\mathcal{E}^{g}_{\theta,\t}(\xi_{\tau})|^4&\leq& c{\rm ess}\sup_{\theta\in \stopo}{\rm ess}\sup_{\tau \in \mathcal{T}_{\theta}}(E[|\xi_{\t}^n-\xi_{\t}|^2|\mathcal{F}_{\theta}])^2\leq c {\rm ess}\sup_{\theta\in \stopo}|U_{\theta}^n|^2,\nonumber
		\end{eqnarray}
		where $U^n$ is given by $U_t^n=E[esssup_{ \tau\in\stopo}|\xi_{\t}^n-\xi_{\t}|^2|\mathcal{F}_{t}]$.
		The process $(U^n_t)_{[t \in [0,T]}$ is right continuous. Thus
		$${\rm ess}\sup_{ \theta \in \stopo}|U_\theta^n|=\sup_{t \in [0,T]}|U_t^n|.$$
		By using this and Doob's martingale inequality in $L^2$, we obtain:
		\begin{eqnarray}\label{expectation} 
			E\left(	{\rm ess}\sup_{\theta\in \stopo}{\rm ess}\sup_{\tau \in \mathcal{T}_{\theta}}|\mathcal{E}^{g}_{\theta,\t}(\xi_{\t}^n)-\mathcal{E}^{g}_{\theta,\t}(\xi_{\tau})|^4\right)\leq cE({\rm ess}\sup_{  \tau\in\stopo}|\xi_{\t}^n-\xi_{\t}|^2).
		\end{eqnarray}
		
		By combining the inequalities (\ref{general-est}) and (\ref{expectation}) with $\|\xi^n-\xi\|_{\mathcal{S}^2} \rightarrow 0$, we derive the desired convergence result.
		\fproof
		\begin{theorem}\label{thm:main}
			Let $g$ be a Lipschitz driver.
			Let $\xi$ be an optional process which  in $\mathcal{S}^2$ which is right upper semicontinuous, and let $(\xi^n)_{n\in\n}$ be given according to (\ref{approx}). Let $Y = \lim_{n\to\infty} Y^n$ in $L^2(\Omega\times[0,T], \mathcal{F}_T\otimes \mathcal{B}([0,T]), P\otimes dt)$,
			then, for each $S \in \stopo$
			\begin{eqnarray}\label{eq2}	
				Y_S={\rm ess}\sup_{\tau\in\mathcal{T}_S}\mathcal{E}^{g}_{S,\tau}(\xi_\t).
			\end{eqnarray} 
		\end{theorem}
		\dproof
		Let $(\xi^n)_{n\in\n}$ be given according to \ref{approx}. From a result of \cite{Agnquenez},  for each stopping time $S$, we have $Y_
		S^n=\mathcal{R}^{g}[\xi^n](S)$. By letting $n$ tend to $\infty$, and using that  the sequence $(Y^n)_{n\in\n}$ converges in $L^2(\Omega\times[0,T], \mathcal{F}_T\otimes \mathcal{B}([0,T]), P\otimes dt)$ a process $Y$ together with Theorem \ref{cornonliear}, we obtain that:
		,
		$$Y_S={\rm ess}\sup_{\tau\in\mathcal{T}_S}\mathcal{E}^{g}_{S,\tau}(\xi_\t)=\mathcal{R}^{g}[\xi](S) $$
		\fproof

		\begin{theorem}\label{princ1}
			Le $g$ be a Lipschitz driver. Let $\xi$ an optional right upper semicontinuous which can be approximated by RCLL barriers $ \xi^n$ monotonically from above in $\mathcal{S}^2$. The reflected BSDE with one reflecting barrier associated with $(g,\xi)$  has a unique solution $(Y,Z,l,A)$. Where  $Y$ is given according to  Theorem \ref{thm:main}.
			$$Y_S={\rm ess}\sup_{\tau\in\mathcal{T}_S}\mathcal{E}^{g}_{S,\tau}(\xi_\t)=\mathcal{R}^{g}[\xi](S)$$
		\end{theorem}
		\dproof
		To prove that $Y$ is the first component of the solution of RBSDE$(\xi,g)$, we apply Theorem \ref{thm:main} combined with Theorem $10.1$ in \cite{grig2}.
		\fproof
		
		In  what follows, the process which aggregates the family $(\mathcal{R}^{g}[\xi](S), S\in \mathcal{T}_0)$ will also be denoted  by $\mathcal{R}^{g}[\xi]$.
		\begin{proposition}\label{snell}
			Let $X$ be an optional process in $\mathcal{S}^2$ such that $\mathcal{R}^{g}[\xi]+X$ is a strong $\mathcal{E}^g$-supermartingale. Let 
			$\bar\xi:=\mathcal{R}^{g}[\xi+X]-X$, then 
			$\mathcal{R}^{g}[\bar \xi]=\mathcal{R}^{g}[\xi]$ a.s. 
		\end{proposition}
		\dproof
		$\mathcal{R}^{g}[\xi+X]$ is the $\mathcal{E}^g$-Snell enveloppe of $\xi+X$. Thus, it is clear that
		$$\xi\leq \mathcal{R}^{g}[\xi+X]-X \;\;\mbox{a.s.} $$
		By using Remark \ref{remark3}, we obtain
		$$\mathcal{R}^{g}[\xi] \leq \mathcal{R}^{g}[\mathcal{R}^{g}[\xi+X]-X]\;\;\mbox{a.s.}$$
		Which yields the first inequality.
		Now, let us prove the second ineqality. 
		We have clearly \\$\xi+X \leq \mathcal{R}^{g}[\xi]+X \;\;\mbox{a.s.}$  
		But $\mathcal{R}^{g}[\xi]+X$ is an $\mathcal{E}^g-$supermartingale by assumption and $\mathcal{R}^{g}[\xi+X]$ is the smallest $\mathcal{E}^g$-supermartingale which is greather than or equal to $\xi+X$. 
		Hence,   
		$$\mathcal{R}^{g}[\xi+X] \leq \mathcal{R}^{g}[\xi]+X\;\;\mbox{a.s.}$$
		It follows that
		$$\mathcal{R}^{g}[\xi+X]-X \leq \mathcal{R}^{g}[\xi]\;\;\mbox{a.s.}$$
		Which yields the desired result:
		$$\mathcal{R}^{g}(\mathcal{R}^{g}[(\xi+X)-X] \leq \mathcal{R}^{g}[\xi]\;\;\mbox{a.s.}$$
		\fproof
		\begin{remark}\label{rema}
			By the property of the Snell envelope, $\xi \leq \bar \xi$. Moroever, if the process $X$ is  continuous, then, $\Delta \bar \xi \leq 0$. This is due to Theorem \ref{princ} and Remark \ref{jumps-left}.
		\end{remark}
		For optional processes $Y$, $Z$, $l$, we set 
		$$g_{Y,Z,l}(t)=g(t,Y_t,Z_t,l_t)\;_;\; \forall t \in [0,T].$$
		\begin{proposition}\label{prop4}
			Let $X$ and $\bar\xi $ are as in Proposition \ref{snell}. Suppose that $X$ is continuous. If $(Y, Z,l, A)$ is the solution of the reflected BSDE associated with $(\bar\xi,g)$. Then, $(Y, Z,l, A)$ is the solution of the reflected BSDE associated with $(\xi, g)$.
		\end{proposition}
		\dproof
		Let $(\tilde Y, \tilde Z,\tilde l, \tilde A)$ be the solution  of the reflected BSDE associated with $(\xi, g_{Y,Z,l})$. Let us prove that $(\tilde Y, \tilde Z,\tilde l,\tilde A)=(Y, Z,l, A)$. By Theorem $10.1$ in \cite{grig2}, $\tilde Y=\mathcal{R}^{g_{Y,Z,l}}[\xi])$ and $ Y=\mathcal{R}^{g_{Y,Z,l}}[\bar\xi]$.  By Proposition \ref{snell}, $Y=\tilde Y$ a.s.
		Moreover, we get by Remark (\ref{rema})
		\[
		\int_0^T(\tilde Y_{t-}-\limsup_{u\uparrow t}\bar \xi_{u})\,d\tilde A_{t}\leq \int_0^T(\tilde Y_{t^-}-\limsup_{u\uparrow t}\xi_{u})\,d \tilde A_t^r
		+\sum_{0<t\le T}(\tilde Y_{t-}-\limsup_{u\uparrow t}\xi_{u})\Delta\tilde A_{t}=0
		\;\mbox{a.s.}\]
		We have also that
		$\sum_{s<T}(\tilde Y_{s}-\bar\xi_{s})\Delta^{+}\tilde A_{s} \leq \sum_{s<T}(\tilde Y_{s}-\xi_{s})\Delta^{+}\tilde A_{s} $.
		Therefore, by uniqueness of the solution $(\tilde Y, \tilde Z,\tilde l, \tilde A)=( Y, Z, l, A)$. This means that $( Y, Z, l,A)$ is the solution  of the reflected BSDE associated with $(\xi,g_{Y,Z,l})$.
		\fproof
		
		In the following theorem,  we give the analoguous of a result of \cite{klimsiak}, in the case when the obstacle process $\xi$ is not necessarily left limitied,  in the setting where the noise is given by a Brownian motion and an independent Poisson measure.
		\begin{theorem}\label{thmfinal}
			Let $\xi$ be a right upper semicontinuous process in $\mathcal{S}^2$, such that
			$\xi_t< \limsup_{u\uparrow t}\xi_{u}$, for all $t\in (0,T]$,  and let $(Y,Z,l,A)$ be the solution of the $RBSDE(\xi,g)$ from Definition \ref{definition1}, then, $(Y^+,Z,l,A^+)$ is the solution of the reflected BSDE with parameters $(\xi^+,g)$. Moreover, for each $S \in \stopo$
			$$Y_{S^+}={\rm ess}\sup_{\tau\geq S}E\Big(\	\xi_\t^++\int^{\tau}_S
			g(s,Y_s^+,Z_s)\,ds|\mathcal{F}_S\Big).$$
		\end{theorem}
		\dproof
		Since $Y\ge \xi$ up to an evanescent set, then of
		course $Y_+\ge \xi_+$ up to an evanescent set. Therefore it is sufficent to show that
		\[
		SK:=\int_0^T ((Y_{t+})_{-}-\limsup_{u\uparrow t}\xi_{u^+})\, dA_{t+}=0.
		\]
		First, let us remark that under the hypothesis 	$\xi_t< \limsup_{u\uparrow t}\xi_{u}$, we have $\limsup_{u\uparrow t}\xi_{u}\leq \limsup_{u\uparrow t}\xi_{u^+}$. Thus,
		\[
		SK\leq \int_0^T(Y_{t-}-\limsup_{u\uparrow t}\xi_{u})\,dA_{t+}=\int_0^T(Y_{t^-}-\limsup_{u\uparrow t}\xi_{u})\,dA^r_t
		+\sum_{0<t\le T}(Y_{t-}-\limsup_{u\uparrow t}\xi_{u})\Delta A_{t+}
		\]
		The first term on the right-hand side is equal to zero since
		$(Y,Z,l,A)$ is the solution of RBSDE associated with $(\xi,g)$.
		Now,  let us prove that the second term is null. We have 
		\[\sum_{0<t< T}(Y_{t-}-\limsup_{u\uparrow t}\xi_{u})\Delta A_{t+}=\sum_{0<t< T}(Y_{t-}-\limsup_{u\uparrow t}\xi_{u})1_{\{Y_t=Y_{t^-}\}}\Delta^+ A_{t}.\]   Suppose that $\Delta^+A_t>0$.
		Then, $Y_t=\xi_t$ by the Skorokhod condition (iv). This together with assumption,
		$\limsup_{u\uparrow t}\xi_{u}> \xi_t$ yields that
		Which completes the proof.
		\fproof
		\begin{corollary}
			Let $Y$ be the first component of the solution of \textnormal{RBSDE}$(\xi,g)$  as in  precedent theorem. Then,
			For each $S \in \stopo$
			$${\rm ess}\sup_{\tau\geq S}E\Big(\	\xi_\t^++\int^{\tau}_S
			g(s,Y_s^+,Z_s,l_s)\,ds|\mathcal{F}_S\Big)={\rm ess}\sup_{\tau> S}E\Big(\	\xi_\t+\int^{\tau}_S
			g(s,Y_s,Z_s,l_s)\,ds|\mathcal{F}_S\Big).$$
		\end{corollary}
		\dproof
		Let $S\in \stopo$, 	Let  us denote  $Y(S)$ the random variable defined by: $$ Y(S)={\rm ess}\sup_{\tau\geq S}E\left[ \xi_\tau+ \int_S^\tau g(u,Y_u,Z_u,l_u) du \mid \FC_{S}\right].$$
		By the same arguments of the proof of Theorem \ref{princ}, the solution $(Y_{t})_{t \in[0,T]}$ of $RBSDE(\xi,g)$ aggregates the family  $( Y(S),\;S\in \stopo)$. Thus, $(Y_{t^+})_{t \in[0,T]}$ aggregates the family $(Y(S^+),\;\;S\in \stopo)$. Now, let 
		$${}^+ Y(S)={\rm ess}\sup_{\tau>S}E\left[ \xi_\tau+ \int_S^\tau g(u,Y_u,Z_u,l_u) du \mid \FC_{S}\right].$$
		Moreover, thanks to a result from optimal stopping theory (cf.  \cite[Proposition 4.14]{kob}) , ${}^+ Y(S)=Y(S^+)$ a.s. Thus, the process $(Y_{t^+})_{t \in[0,T]}$ aggregates the family $({}^+ Y(S),\;S\in \stopo)$.
		Hence, 
		$Y_{S^+}={\rm ess}\sup_{\tau> S}E\Big(\xi_\t+\int^{\tau}_S
		g(s,Y_s,Z_s)\,ds|\mathcal{F}_S\Big)$. By Theorem \ref{thmfinal},$$Y_{S^+}={\rm ess}\sup_{\tau\geq S}E\Big(\	\xi_\t^++\int^{\tau}_S
		g(s,Y_s^+,Z_s^+)\,ds|\mathcal{F}_S\Big).$$
		This yields the desired result.
		\fproof
		\begin{remark}
			In particular, if $\xi$ is a right upper semicontinuous optional process satisfying $\xi_t \leq \limsup_{u\uparrow t}\xi_t$, then 
			$${\rm ess}\sup_{\tau\geq S}E\left(\xi_\t^+|\mathcal{F}_S\right)={\rm ess}\sup_{\tau> S}E\left(\xi_\t|\mathcal{F}_S\right).$$
		\end{remark}
		\begin{lemma}
			Let $\bar \xi$ is given as in Proposition \ref{prop4}. Let  $(Y^n,Z^n,A^n)$ be the solution of the following BSDE:
			$$Y^n_t=\bar \xi _T+\int_t^T g(u, Y_u^n, Z_u^n)du+\int_t^Tn(Y_u^n-\bar \xi_{u}^+)^-ds - \int_t^T Z_u^n d W_u -\int_{t}^T\int_E  l_u^n(e) \tilde N(du,de) .$$
			Let $$\bar Y^n:=\xi\vee Y_t^n. $$ Then, $$\bar Y_t^n \uparrow Y_t,\;\;t \in [0,T],$$
			where  $Y$ is the first component of the solution  $(Y,Z,l,A)$ of $RBSDE(\xi,g)$.			
		\end{lemma}
		\dproof
		Let $(\bar Y,\bar Z, \bar l,\bar A)$ be the solution of $RBSDE(\bar\xi,g)$. By Proposition \ref{prop4} $(\bar Y,\bar Z,\bar l,\bar A)=(Y,Z,l,A)$. We have by Remark \ref{rema} that $\Delta \bar \xi \leq 0$.  Then, by Theorem \ref{thmfinal},  $(Y^+,Z,l,A^+)$ is the solution of $RBSDE(\bar \xi^+,g)$. On the other hand, by Hamad\`ene and Ouknine \cite{hamadeneouknine03},  $ Y^n \uparrow  Y^+ .$ Hence, $\bar Y_t^n\uparrow  \xi_t \vee Y_t^+,\;\;t \in [0,T]$. The result follows from Proposition \ref{prop-fon}.
		\fproof
		\section{Generalized Skorohod 
			condition in  the sens of Peng-Xu}
		To show how a solution of BSDE can be reflected by a very irregular $L^2$ obstacle, Peng and Xu \cite{Peng2}, found a new fromulation of Skorokhod condition. 
		\begin{definition}[RBSDEs in the sens of Peng-Xu]
			Let $g
			$ be a driver, $\xi$ an obstacle.We say a triple of processes $(Y, Z, A)$ is  a solution of the reflected BSDE with standard parameters $(g,\xi)$ if
			$(Y, Z, A)\in \mathbb{D}^2\times L^2\times\mathbb{D}^2$,
			$$Y_\tau = \xi_T + \int_\tau^T g(\cdot, t, Y_t, Z_t) dt - \int_\tau^T Z_t d W_t + A_T - A_\tau\quad\mbox{a.s. for all}\quad \tau\in\mathcal{T}.$$
			$Y\ge \xi$ $\;\;\;\;\;dt \otimes dP$,\\
			the following generalized Skorohod condition holds
			$$\int_0^T(Y_t-\xi_t^*) d A_t = 0\quad\mbox{a.s. for all}\quad \xi^* \in \mathbb{D}^2\quad \mbox{such that} \quad\xi_t \leq \xi_t^* \leq Y_t \quad 
			a.s., a.e.$$	
		\end{definition} 
		In \cite{Peng2}, Peng and Xu shows the existence of a solution of a unique solution $(Y,Z,A)$ by penalization method. The aim of this part is to give a new appraoch which avoids the generalized Skorohd condition involving $\xi^*$.	
		\begin{theorem}\label{rexiuni}
			Let $\xi$ be an optional process in $\mathcal{S}^2$ such that
			$\xi\le \xi^+$ up to an evanescent set,  and let $(Y,Z,A)$ be the solution of the reflected BSDE with parameters $(\xi,g)$ from Definition \ref{definition1}, then, $(Y,Z,A)$ is the solution of the reflected BSDE with parameters $(\xi,g)$ in the sens of Peng-Xu.
		\end{theorem}
		
		\dproof
		First, note that the process $Y$ is right continuous. This follows from the assumption $\xi\le \xi^+$ up to an evanescent set and an application of Lemma \ref{lem right}. Let $\xi^*$ be a cadlag process, such that $\xi \leq \xi^* \leq Y$  $dt\otimes dP$. We need show that
		\[
		SK:=\int_0^T ((Y_{t+})_{-}-\xi_{t^-}^*)\, dA_t=0.
		\]
		Since $\xi^* \geq  \xi$ $dt\otimes dP$ , we get 
		\begin{eqnarray}
			SK=\int_0^T(Y_{t^-}- \xi_{t^-}^*)dA_t\leq \int_0^T(Y_{t^-}-\limsup_{u\uparrow t}\xi_{u})dA_t.
		\end{eqnarray}	
		The term on the right hand side is equal to zero since $(Y,Z,A)$ is the solution of $RBSDE(\xi,g)$. Hence the result.
		\fproof

	\end{document}